\documentclass[12pt,reqno,a4wide]{amsart}

\usepackage{tikz} 
\usepackage{calrsfs}
\usepackage{mathrsfs}

\usepackage{array}
\usepackage{hyperref}

\usetikzlibrary{decorations.pathreplacing}
\usetikzlibrary{fit}

\usepackage{pgfplots}

\allowdisplaybreaks

\newtheorem{theorem}{Theorem}

\catcode`,\active

\catcode`\,12

\begin{document}
\bibliographystyle{plain}

%
%

	\title
	{Peakless Motzkin paths of bounded height}

	\author[H. Prodinger ]{Helmut Prodinger }
	\address{Department of Mathematics, University of Stellenbosch 7602, Stellenbosch, South Africa
	and
NITheCS (National Institute for
Theoretical and Computational Sciences), South Africa.}
	\email{hproding@sun.ac.za}

	
	\begin{abstract}
		There was recent interest in Motzkin paths without peaks (peak: up-step followed immediately by down-step); additional results 
		about this interesting family is worked out. The new results are the enumeration of such paths that live in a strip $[0..\ell]$, and as
		consequence the asymptotics of the average height, which is given by $2\cdot 5^{-1/4}\sqrt{\pi n}$. Methods include the kernel method
		and singularity analysis of generating functions.
	\end{abstract}
	

\maketitle

\section{Introduction}

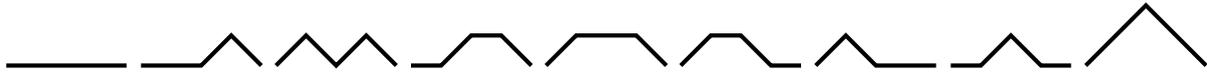
\begin{figure}[h]
	\label{all1}
\begin{center}

			\begin{tikzpicture}[scale=0.4]
				
				\draw[ultra thick] (0,0) to (1,0) to (2,0) to (3,0) to (4,0)  ;
			\end{tikzpicture}
			\begin{tikzpicture}[scale=0.4]
				
				\draw[ ultra thick] (0,0) to (1,0) to (2,0) to (3,1) to (4,0)  ;
			\end{tikzpicture}
			\begin{tikzpicture}[scale=0.4]
				
				\draw[ ultra thick] (0,0) to (1,1) to (2,0) to (3,1) to (4,0)  ;
			\end{tikzpicture}
			\begin{tikzpicture}[scale=0.4]
				
				\draw[ ultra thick] (0,0) to (1,0) to (2,1) to (3,1) to (4,0)  ;
			\end{tikzpicture}
			\begin{tikzpicture}[scale=0.4]
				
				\draw[ ultra thick] (0,0) to (1,1) to (2,1) to (3,1) to (4,0)  ;
			\end{tikzpicture}
			\begin{tikzpicture}[scale=0.4]
				
				\draw[ ultra thick] (0,0) to (1,1) to (2,1) to (3,0) to (4,0)  ;
			\end{tikzpicture}
			\begin{tikzpicture}[scale=0.4]
				
				\draw[ ultra thick] (0,0) to (1,1) to (2,0) to (3,0) to (4,0)  ;
			\end{tikzpicture}
			\begin{tikzpicture}[scale=0.4]
				
				\draw[ ultra thick] (0,0) to (1,0) to (2,1) to (3,0) to (4,0)  ;
			\end{tikzpicture}
			\begin{tikzpicture}[scale=0.4]
				
				\draw[ ultra thick] (0,0) to (1,1) to (2,2) to (3,1) to (4,0)  ;
			\end{tikzpicture}
	 	
\end{center}

\caption{All 9 Motzkin of 4 steps (length 4).}
\end{figure}

Motzkin paths are cousins of the more famous Dyck paths. They appear first in \cite{Motzkin48}. In the encyclopedia \cite{OEIS} they 
are sequence A001006, with many references given. They consist of up-steps $U=(1,1)$, down-steps $D=(1,-1)$ and horizontal (flat)
steps $F=(1,0)$. Slightly different notations are also in use. They start at the origin and must never go below the $x$-axis. Usually one
requires the path to end on the $x$-axis as well, but occasionally one  uses the term \emph{Motzkin path} also for paths that end
on a different level.  Figure~\ref{all1} shows all Motzkin paths of 4 steps (=length 4).

An important concept is the \emph{height} of a path. It is the maximal $y$-coordinate  when scanning the path (from left to right, say). For
the paths in Figure~\ref{all1}, the heights are (in this order) $0,1,1,1,1,1,1,1,2$. The average height of all Motzkin paths of length $n$ was computed in an
early paper of the present writer \cite{prodinger-ars}.

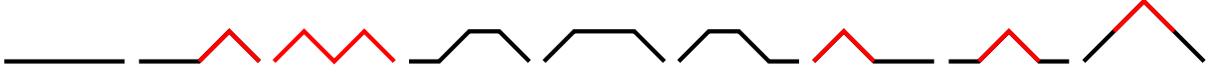
\begin{figure}
	\begin{center}

		\begin{tikzpicture}[scale=0.4]
			
			\draw[ultra thick] (0,0) to (1,0) to (2,0) to (3,0) to (4,0)  ;
		\end{tikzpicture}
		\begin{tikzpicture}[scale=0.4]
			
			\draw[ ultra thick] (0,0) to (1,0) to (2,0) to (3,1) to (4,0)  ;
				\draw[ ultra thick,red]  (2,0) to (3,1) to (4,0)  ;
		\end{tikzpicture}
		\begin{tikzpicture}[scale=0.4]
			
			\draw[ ultra thick,red] (0,0) to (1,1) to (2,0) to (3,1) to (4,0)  ;
		\end{tikzpicture}
		\begin{tikzpicture}[scale=0.4]
			
			\draw[ ultra thick] (0,0) to (1,0) to (2,1) to (3,1) to (4,0)  ;
		\end{tikzpicture}
		\begin{tikzpicture}[scale=0.4]
			
			\draw[ ultra thick] (0,0) to (1,1) to (2,1) to (3,1) to (4,0)  ;
		\end{tikzpicture}
		\begin{tikzpicture}[scale=0.4]
			
			\draw[ ultra thick] (0,0) to (1,1) to (2,1) to (3,0) to (4,0)  ;
		\end{tikzpicture}
		\begin{tikzpicture}[scale=0.4]
			
			\draw[ ultra thick] (0,0) to (1,1) to (2,0) to (3,0) to (4,0)  ;
				\draw[ ultra thick, red] (0,0) to (1,1) to (2,0);
		\end{tikzpicture}
		\begin{tikzpicture}[scale=0.4]
			
			\draw[ ultra thick] (0,0) to (1,0) to (2,1) to (3,0) to (4,0)  ;
						\draw[ ultra thick, red] (1,0) to (2,1) to (3,0);
		\end{tikzpicture}
	 		\begin{tikzpicture}[scale=0.4]
			
			\draw[ ultra thick] (0,0) to (1,1) to (2,2) to (3,1) to (4,0)  ;
			\draw[ ultra thick, red]  (1,1) to (2,2) to (3,1)   ;
		\end{tikzpicture}
	\end{center}

	\caption{Motzkin paths with peaks indicated.}\label{alla}

\end{figure}

Recently, I learned from the paper \cite{wall} that there is interest in \emph{peakless Motzkin paths.} A peak in a Motzkin path is a sequence of an up-step followed
immediately by a down-step. Figure~\ref{alla} indicates all peaks in the list of Motzkin paths of length 4.
The enumerating sequence is A004148 in \cite{OEIS}, where one can find several references; one recent paper about the subject is \cite{cameron-sullivan}.
A general discussion about forbidden patterns, concentrating on analytic aspects, is in \cite{asi-pattern}.

\begin{figure}[h]
	
	\begin{center}

		\begin{tikzpicture}[scale=0.4]
			
			\draw[ultra thick] (0,0) to (1,0) to (2,0) to (3,0) to (4,0)  ;
		\end{tikzpicture}
		\begin{tikzpicture}[scale=0.4]
			
			\draw[ ultra thick] (0,0) to (1,0) to (2,1) to (3,1) to (4,0)  ;
		\end{tikzpicture}
		\begin{tikzpicture}[scale=0.4]
			
			\draw[ ultra thick] (0,0) to (1,1) to (2,1) to (3,1) to (4,0)  ;
		\end{tikzpicture}
		\begin{tikzpicture}[scale=0.4]
			
			\draw[ ultra thick] (0,0) to (1,1) to (2,1) to (3,0) to (4,0)  ;
		\end{tikzpicture}
		
	\end{center}
	\caption{Peakless Motzkin paths of length 4.}	\label{all3}

\end{figure}
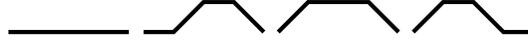

The paths without peaks are called peakless, and there are four of them, as shown in Figure~\ref{all3}.

\section{A warmup: Enumeration of peakless Motzkin paths via the kernel method}

 \begin{figure}[h]

 	\begin{center}
 		\begin{tikzpicture}[scale=1.5,main node/.style={circle,draw,font=\Large\bfseries}]

 			\foreach \x in {0,1,2,3,4,5,6,7,8}
 			{
 				\draw (\x,0) circle (0.05cm);
 				\fill (\x,0) circle (0.05cm);
 				\draw (\x,-1) circle (0.05cm);
 				\fill (\x,-1) circle (0.05cm);
 			}

 			\fill (0,0) circle (0.08cm);

 			\foreach \x in {0,2,4,6}
 			{
 			}
 			\foreach \x in {0,...,8}
 			{
			\draw[thick,-latex ] (\x,0)  ..  controls (\x-0.25,0.5) and  (\x+0.25,0.5) .. (\x,0) ;	
 			}

 			\foreach \x in {0,1,2,3,4,5,6,7}
 			{
 				\draw[thick, latex-] (\x,0) to  (\x+1,0);	
\draw[thick, -latex] (\x,-1) to  (\x+1,-1);	
 				\draw[thick,red, -latex] (\x,0) to  (\x+1,-1);	
 				\node at  (\x+0.2,0.15){\tiny$\x$};
 			}			
 			
 			\foreach \x in {0,1,2,3,4,5,6,7,8}
 			{
 				\draw[thick, -latex] (\x,-1) to  (\x,0);	
 			}

 			\node at  (8+0.2,0.15){\tiny$8$};

 		\end{tikzpicture}
 	\end{center}
 	\caption{Graph (automaton) to recognize peakless Motzkin paths. Starting  at the origin and ending at nodes labelled 0 corresponds to Motzkin paths, and
 	ending at a node labelled $k$ to  a path that ends at level $k$. }
 	\label{purpel}
 \end{figure}
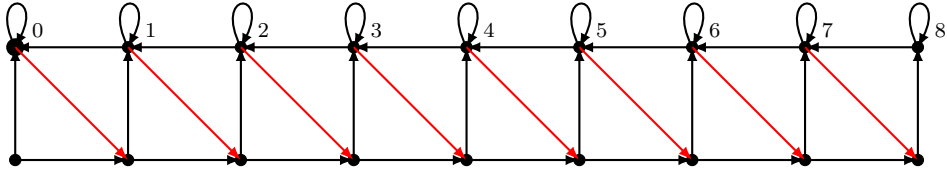

We will use the following generating functions: $[z^n]f_i(z)$ is the number of peakless paths ending at state $i$ in the top layer (Figure~\ref{purpel}),
$[z^n]g_i(z)$ is the number of peakless paths ending at state $i$ in the bottom layer; for convenience, we mostly write $f_i$ and $g_i$.
The following recursions can be read off the automaton, by considering the last step separately.
 \begin{align*}
f_0&=1+zf_0+zf_1+zg_0,\\
f_i&=zf_{i}+zf_{i+1}+zg_i,\quad i\ge1,\\
g_0&=0,\\
g_{i+1}&=zf_i+zg_i,\quad i\ge0.
 \end{align*}
To solve this system, one introduces double generating functions
\begin{equation*}
F(u,z)=\sum_{i\ge0}u^if_i(z)\quad\text{and}\quad G(u,z)=\sum_{i\ge0}u^ig_i(z).
\end{equation*}
Again, for convenience, we mostly write $F(u)$ and $G(u)$. By summing the recursions, we find
\begin{align*}
F(u)&=1+zF(u)+\frac zu(F(u)-F(0))+zG(u),\\
G(u)&=zuF(u)+zuG(u)=\frac{zuF(u)}{1-zu}.
\end{align*}
Eliminating one function, we are left to solve (note that $F(0)=f_0$)
\begin{equation*}
	F(u)=1+zF(u)+\frac zu(F(u)-F(0))+\frac{z^2uF(u)}{1-zu}.
\end{equation*}
Rewriting the functional equation, we find
\begin{equation*}
F(u)=\frac{(-u+zF(0))(1-zu)}{zu^2+(z-z^2-1)u+z}=\frac{(-u+zF(0))(1-zu)}{z(u-s_1)(u-s_2)}
\end{equation*}
with
\begin{equation}\label {ss1}
s_1=\frac{1-z+z^2+\sqrt{(1+z+z^2)(1-3z+z^2)}}{2z}
\end{equation}
and
\begin{equation}\label {ss2}
	s_2=\frac{1-z+z^2-\sqrt{(1+z+z^2)(1-3z+z^2)}}{2z}.
\end{equation}
Note that $s_1s_2=1$. Plugging $u=0$ into that equation does not help, but one of the factors from the denominator can be cancelled.
This is (a simple instance of) the kernel method. Some twenty years ago, I collected various related examples \cite{prodinger-kernel}, and many more just recently~\cite{prodinger-garden}.

Since $s_2=z+z^2+\cdots$, the factor $(u-s_2)$ must cancel, since $1/(u-s_2)$ would not have a power series expansion around $u$, $z$ close to zero.
Performing the cancellation, we get
\begin{equation*}
F(u)=\frac{zs_2-1+zu-z^2F(0)}{z(u-s_1)}\quad\text{and}\quad
F(0)=\frac{zs_2-1-z^2F(0)}{-zs_1}.
\end{equation*}
Solving,
\begin{align*}
F(0)=f_0&=\frac{s_2}{z}=\frac{1-z+z^2-\sqrt{(1+z+z^2)(1-3z+z^2)}}{2z^2}\\&=1+z+z^2+2 z^3+4 z^4+8 z^5+17 z^6+\cdots
\end{align*}
We are mostly interested in
\begin{equation*}
H(u):=F(u)+G(u)=
\frac {  -u+zF(0) }{zu+z-{z}^{2}u-u+z{u}^{2}},
\end{equation*}
although $F(u)$ and $G(u)$ could be computed separately as well. The coefficients of $[u^k]H(u)$ are just enumerating all peakless Motzkin paths ending on level $k$. Hence
\begin{equation*}
H(u)=\frac{-u+zF(0)}{zu+z-{z}^{2}u-u+z{u}^{2}}=\frac{-1}{z(u-s_1)}.
\end{equation*}
Expanding and noting that $s_1=1/s_2$, we find
\begin{equation*}
[u^k]H(u)=\frac{s_2^{k+1}}{z}.
\end{equation*}
This could be seen as well by a canonical decomposition of a peakless Motzkin path according to the last return to the $x$-axis.

The denominator of $H(u)$ has a special significance; if we write $h_k=[u^k]H(u)$, the recursion for these quantities can be read off from the denominator:
\begin{equation}\label{nenn}
zh_k+(z-z^2-1)h_{k-1}+zh_{k-2}=0,
\end{equation}
this can be checked directly as well by inserting $h_k=s_2^{k+1}/z$ and simplifying.

The sequence enumerating peakless Motzkin paths is A004148 in \cite{OEIS}. If we call them $m(n)=[z^n]s_2/z$, then the software \textsf{Gfun}, implemented in Maple,
produces the recursion
\begin{equation*}
	nm  ( n  ) -  ( 2n+3  ) m  ( n+1 ) -  (n+3  ) m  ( n+2  ) -  (2n+9 ) m  ( n+3  ) +  ( n+6  ) m  ( n+4 )=0,
\end{equation*}
with initial values $m  ( 0  ) =1$, $m  ( 1  ) =1$, $m  ( 2) =1$, $m  ( 3  ) =2  $.

A second order linear recursion with constant coefficients is driven by the characteristic equation and its two roots. Not surprisingly, they  are $s_1$ and $s_2$.

For completeness, we mention the asymptotics of the coefficients $m(n)$ of
\begin{equation*}
\frac{s_2}{z}=\frac{1-z+z^2-\sqrt{(1+z+z^2)(1-3z+z^2)}}{2z^2}.
\end{equation*}
This is a standard application of \emph{singularity analysis of generating functions}, as described in \cite{FS}. First, we must consider the closest singularity to the origin.
The candidates are the solutions of $(1+z+z^2)(1-3z+z^2)=0$. There are two complex solutions of absolute value 1, which are irrelevant, and then $\frac{3\pm\sqrt5 }{2}$.
The relevant value is $\varrho=\frac{3-\sqrt5 }{2}$; note that $1/\varrho=\frac{3+\sqrt5 }{2}$, which is the square of the golden ratio $\phi=\frac{1+\sqrt5}{2}$ from the Fibonacci fame.

The local expansion around $z\sim\varrho$ looks like
\begin{equation*}
\frac{s_2}{z}\sim \frac1{\varrho}- \frac{5^{1/4}}{\varrho}\sqrt{1-\frac{z}{\varrho}},
\end{equation*}
and following the principles of singularity analysis we might translate this to the coefficients:
\begin{equation*}
[z^n]\frac{s_2}{z}\sim \frac{5^{1/4}\varrho^{-n-1}}{2\sqrt\pi n^{3/2}}.
\end{equation*}
\section{Peakless Motzkin paths of bounded height}

We fix a parameter $\ell\ge0$ and postulate that $[u^j]H(u)=0$ for $j>\ell$. This means that states $\ell+1,\ell+2,\dots$ (on both layers) can never be reached.
The recursion, compare (\ref{nenn})
\begin{equation*}
zh_{k+1}+(z-z^2-1)h_{k}+zh_{k-1}=0
\end{equation*}
is then best written as a matrix equation:
\begin{equation*}
\begin{pmatrix} z-z^2-1& z&0&0&\dots\\ 
	z&z-z^2-1& z&0&0&\dots\\ 
0&	z&z-z^2-1& z&0&\dots\\ 
\vdots&\vdots&\vdots&\ddots&\ddots&\\
&&& &z&z-z^2-1
	\end{pmatrix}
\begin{pmatrix} h_0\\h_1\\ h_2\\ \vdots \\    h_\ell
\end{pmatrix}
=\begin{pmatrix} -1\\0\\ 0\\ \vdots \\    0
\end{pmatrix}
\end{equation*}

Let $\mathcal{D}_\ell$ be the determinant of the $(\ell+1)\times(\ell+1)$ matrix. $\mathcal{D}_0=z-z^2-1$, $\mathcal{D}_1=(z-z^2-1)^2-z^2=(1-z)^2(1+z^2)$.
It is a bit easier to work with $\mathcal{D}_{-1}=1$ instead of $\mathcal{D}_{1}$. The solution is, by standard methods,
\begin{equation*}
\mathcal{D}_\ell=-\frac1W(-zs_1)^{\ell+2}+\frac1W(-zs_2)^{\ell+2},
\end{equation*}
where we use the abbreviation $W=\sqrt{(1+z+z^2)(1-3z+z^2)}$.

The quantity $h_0$, describing \emph{all} paths, restricted as described, can, by Cramer's rule, be written as
\begin{align*}
h_0=\frac{-\mathcal{D}_{\ell-1}}{\mathcal{D}_\ell}&=-\frac{-\frac1W(-zs_1)^{\ell+1}+\frac1W(-zs_2)^{\ell+1}}{-\frac1W(-zs_1)^{\ell+2}+\frac1W(-zs_2)^{\ell+2}}\\
	&=\frac1z\frac{-(-s_1)^{\ell+1}+(-s_2)^{\ell+1}}{(-s_1)^{\ell+2}-(-s_2)^{\ell+2}}
=\frac1z\frac{s_1^{\ell+1}-s_2^{\ell+1}}{s_1^{\ell+2}-s_2^{\ell+2}}.
\end{align*}
In the limit $\ell\to\infty$, $h_0=\dfrac1{zs_1}=\dfrac{s_2}{z}$. This is, as we have seen already, the generating function of peakless Motzkin paths without boundary. The other functions $h_i$ could be computed by Cramer's rule as well, but we concentrate only on
the paths that return to the origin and are bounded by $\ell$. They live in the strip $[0..\ell]$. At this stage, we drop the `$0$' from the notation and make the `$\ell$' explicit by writing $A_{n,\ell}$. We summarize the results:
\begin{theorem}
	The number of peakless Motzkin paths of length $n$ (returning to the $x$-axis), bounded by $\ell$, is given by
	\begin{equation*}
A_{n,\ell}=[z^n]\frac1z\frac{s_1^{\ell+1}-s_2^{\ell+1}}{s_1^{\ell+2}-s_2^{\ell+2}},
	\end{equation*}
with the functions $s_1$ and $s_2 $ given in (\ref{ss1}) and (\ref{ss2}). The formula is only correct for $\ell\ge1$; if $\ell=0$, there is only one such path
of length of $n$, namely consisting of flat steps only; the generating function must then be replaced by $\dfrac{1}{1-z}$.
\end{theorem}
We will use the notations ($\ell\ge1$)
\begin{equation*}
	A_{\ell}(z)=\frac1z\frac{s_1^{\ell+1}-s_2^{\ell+1}}{s_1^{\ell+2}-s_2^{\ell+2}}\quad\text{and}\quad A_{\infty}(z)=\dfrac{s_2}{z}.
\end{equation*}
As can be checked directly (best with a computer), there is the recursion of the continued fraction type
\begin{equation*}
A_{\ell}=\frac1{1-z+z^2-z^2A_{\ell-1}}, 
\end{equation*}
and thus
\begin{equation*}
A_1=\cfrac{1}{1-z+z^2-\cfrac{z^2}{1-z+z^2}},\quad A_2=\cfrac{1}{1-z+z^2-\cfrac{z^2}{1-z+z^2-\cfrac{z^2}{1-z+z^2}}}
\end{equation*}
and so on. As discussed, the formula for $A_0=\frac1{1-z}$ is slightly different.

Continued fraction expansions are very common in the context of generating functions of lattice paths of bounded height; the first example is (perhaps) \cite{BrKnRi72}.

The total generating function $A_\infty$ has a pretty expansion as well, viz.
\begin{equation*}
A_\infty=1+\cfrac{z}{1-\cfrac{z}{1-\cfrac{z}{1-\cfrac{z^3}{1-\cfrac{z}{1-\cfrac{z}{1-\cfrac{z^3}{1-\ddots}}}}}}}	
\end{equation*}
this can be checked directly. I am wondering whether this might have an easy combinatorial interpretation.

\medskip

Now we consider the \emph{average height} of peakless Motzkin paths of length $n$, assuming that all of them are equally likely. As always, when enumerating the average height, the relevant formula is
\begin{equation*}
\frac{[z^n]\sum_{\ell\ge0}\big(A_\infty(z)-A_{\ell}(z)\big)}{[z^n]A_\infty(z)}.
\end{equation*}

The difference $A_{\infty}(z)-A_{\ell}(z)$ enumerates paths of height $>\ell$, or
\begin{equation*}
\dfrac{s_2}{z}-\frac1z\frac{s_1^{\ell+1}-s_2^{\ell+1}}{s_1^{\ell+2}-s_2^{\ell+2}}=\frac1z\frac{(1-s_2^{2})s_2^{\ell+1}}{s_1^{\ell+2}-s_2^{\ell+2}}
=\frac{W}{z^2}\frac{s_2^{\ell+2}}{s_1^{\ell+2}-s_2^{\ell+2}}
\end{equation*}
This must be expanded around $z=\varrho$; we decrease $\ell$ by one, since then we enumerate paths of height $\ge\ell$. Since $s_2\sim1$, we find the approximation
\begin{equation*}
\frac{W}{\varrho^2}\frac{1}{s_1^{\ell+1}-1}=\frac{W}{\varrho^2}\frac{s_2^{\ell+1}}{1-s_2^{\ell+1}}=\frac{W}{\varrho^2}\sum_{k\ge1}s_2^{(\ell+1)k}
\end{equation*}
A local expansion yields
\begin{equation*}
\frac{W}{\varrho^2}\sim\frac{2\cdot5^{1/4}}{\varrho}\Big(1-\frac z\varrho\Big)^{1/2}.
\end{equation*}
From our earlier computation,
\begin{equation*}
s_2\sim 1 -5^{1/4}\Big(1-\frac z\varrho\Big)^{1/2}.
\end{equation*}
 From \cite{HPW} we conclude that
\begin{equation*}
\sum_{k\ge1}\frac{s_2^k}{1-s_2^k}\sim -\frac{\log(1-s_2)}{1-s_2},
\end{equation*}
as $s_2\to 1$. Our paper \cite{HPW} has many more technical details about a similar scenario. Putting both expansions together (we don't care about a missing term in the sum as we only want to work out the leading term),
\begin{align*}
\frac{W}{\varrho^2}\sum_{k,\ell\ge1}s_2^{(\ell+1)k}&\sim \frac{2\cdot5^{1/4}}{\varrho}\Big(1-\frac z\varrho\Big)^{1/2}\cdot\frac{-\log(1-s_2)}{1-s_2}\\
&\sim\frac{2\cdot5^{1/4}}{\varrho}\Big(1-\frac z\varrho\Big)^{1/2}\cdot\frac{-\log\Big(5^{1/4}\big(1-\dfrac z\varrho\big)^{1/2}\Big)}{5^{1/4}\big(1-\dfrac z\varrho\big)^{1/2}}\\
&\sim-\frac{2}{\varrho}\log\Big(5^{1/4}\big(1-\dfrac z\varrho\big)^{1/2}\Big)
\sim-\frac{2}{\varrho}\log\Big(1-\dfrac z\varrho\Big)^{1/2}
\sim-\frac{1}{\varrho}\log\Big(1-\dfrac z\varrho\Big).
\end{align*}
By singularity analysis (transfer theorem) we find that
\begin{equation*}
[z^n]\frac{W}{\varrho^2}\sum_{k,\ell\ge1}s_2^{(\ell+1)k}\sim-[z^n]\frac{1}{\varrho}\log\Big(1-\dfrac z\varrho\Big)\sim \frac{\varrho^{-n-1}}{n}.
\end{equation*}
As discussed before, the total number of peakless Motzkin paths of length $n$ is asymptotic to
\begin{equation*}
	[z^n]\frac{s_2}{z}\sim \frac{5^{1/4}\varrho^{-n-1}}{2\sqrt\pi n^{3/2}}.
\end{equation*}
and for the average height we have to consider the quotient of the last two expressions, which is 
\begin{equation*}
	\frac{\varrho^{-n-1}}{n}\frac{2\sqrt\pi n^{3/2}}{5^{1/4}\varrho^{-n-1}}	=\frac{2\sqrt{\pi n} }{5^{1/4}}.
\end{equation*}
The numerical constant $2/5^{1/4}=1.337480610$. This can be compared with the average height of \emph{all} Motzkin paths of length $n$, which is
asymptotic to $\sqrt{\frac{\pi n}3}$, see \cite{Prodinger-three}; $3^{-1/2}=0.5773502693$.

\bibliographystyle{plain}

\end{document}